\documentclass[11pt, twoside, a4paper, fleqn]{article}

\usepackage{latexsym}
\usepackage{amscd}
\usepackage{theorem}
\usepackage{pifont}
\usepackage{amsfonts}
\usepackage{xspace}
\usepackage{amssymb}
\usepackage{fancyhdr}
\usepackage{mathrsfs}
\usepackage{amsmath}%
\usepackage{graphics}%
\usepackage{graphicx}%
\usepackage{euscript}%
\usepackage{psfrag}%
\usepackage{empheq}%
\usepackage{upgreek}
\usepackage[dvipsnames]{xcolor}

\DeclareMathAlphabet{\mathpzc}{OT1}{pzc}{m}{it}

{\theorembodyfont{\slshape} \newtheorem{theorem}{Theorem}[section]}
{\theorembodyfont{\slshape} \newtheorem{definition}[theorem]{Definition}}
{\theorembodyfont{\slshape} }
{\theorembodyfont{\slshape} }
{\theorembodyfont{\slshape} \newtheorem{corollary}[theorem]{Corollary}}
{\theorembodyfont{\slshape} \newtheorem{remark}[theorem]{Remark}}
{\theorembodyfont{\slshape} \newtheorem{example}[theorem]{Example}}

\newcommand{\complex}{\mathbb{C}}

\newenvironment{proof}{\noindent\textbf{Proof:\ }}{$\hfill{\Box}$}

\numberwithin{equation}{section}

\title{\textsc{Dynamics of Products of Nonnegative Matrices}}    

\author{Sachindranath Jayaraman \\ {\tt sachindranathj@iisertvm.ac.in;\ sachindranathj@gmail.com} \bigskip 
	\\ Yogesh Kumar Prajapaty \\ {\tt prajapaty0916@iisertvm.ac.in} \bigskip 
	\\ Shrihari Sridharan \\ {\tt shrihari@iisertvm.ac.in} \bigskip 
 \bigskip \\ {\sl Indian Institute of Science Education and Research}
  \\ {\sl Thiruvananthapuram (IISER-TVM), India.}} 

\date{\today}

\begin{document}

\maketitle
\thispagestyle{empty}

\begin{abstract} 
\noindent 
The aim of this manuscript is to understand the dynamics of products of nonnegative matrices. We extend a well known consequence of the Perron-Frobenius theorem on the periodic points of a nonnegative matrix to products of finitely many nonnegative matrices associated to a word and later to products of nonnegative matrices associated to a word, possibly of infinite length. We also make use of an appropriate definition of the exponential map and the logarithm map on the positive orthant of $\mathbb{R}^{n}$ and explore the relationship between the periodic points of certain subhomogeneous maps defined through the above functions and the periodic points of matrix products, mentioned above. 
\end{abstract} 

\begin{tabular}{r c l} 
{\bf Keywords} & : & Products of nonnegative matrices; \\ 
& & Common eigenvectors and common periodic points; \\ 
& & Subhomogeneous maps; \\ 
& & Intersecting orbits of infinite matrix products. \\ 
& & \\ 
& & \\ 
{\bf MSC 2010 Subject} & : & 15B48; 15A27; 37H12. \\ 
{\bf Classifications} & & 
\end{tabular} 
\bigskip 

\section{Introduction} 

\noindent 
Given a collection $\mathcal{F}$ of functions on a set $\Omega$, an element $w \in \Omega$ is said to be a common fixed point for $\mathcal{F}$ if $f(w) = w$ for all $f \in \mathcal{F}$. The existence and computation of such a point has been a topic of interest among several mathematicians. Of particular interest is when the collection is a multiplicative semigroup or a group $\mathcal{M}$ of matrices, where a more general question on the existence of common eigenvectors arises. A classic example of a multiplicative semigroup of matrices is the collection of matrices whose entries are nonnegative real numbers. In a recent work, Bernik {\it et al} \cite{bdklmmor} determined certain conditions that ensures the existence of a common fixed point and more generally the existence of a common eigenvector for such a collection $\mathcal{M}$. The existence of common eigenvectors for a collection of matrices is in itself a nontrivial question and plays a major role in many problems in matrix analysis. For recent results on periods and periodic points of iterations of sub-homogeneous maps on a proper polyhedral cone, we refer the reader to \cite{agln} (for instance, see Theorem (4.2)) and the references cited therein.
\medskip 

\noindent 
We work throughout with the field $\mathbb{R}$ of real numbers. Let $M_{n} (\mathbb{R})$ denote the real vector space of $n \times n$ matrices. The subset of $M_{n} (\mathbb{R})$ consisting of matrices whose entries are nonnegative real numbers (such a matrix is usually called a \emph{nonnegative matrix}) is denoted by $M_{n} (\mathbb{R}_{+})$. For any matrix $A \in M_{n} (\mathbb{R})$, we denote and define \emph{the spectrum, the spectral radius} and \emph{the norm} of $A$ respectively, as follows: 
\begin{eqnarray*} 
{\rm spec} (A) & = & \text{the set of all eigenvalues of $A$, some of which may be complex numbers}; \\ 
\rho (A) & = & \max \big\{ | \lambda |\ :\ \lambda \in {\rm spec} (A) \big\}; \\ 
\Vert A \Vert  & = & \text{the operator norm of $A$, induced by the Euclidean norm of $\mathbb{R}^{n}$}. 
\end{eqnarray*} 
For any $N < \infty$, we fix a finite collection of matrices, $\big\{ A_{1},\, A_{2},\, \cdots,\, A_{N}\ :\ A_{r} \in M_{n} (\mathbb{R}_{+}) \big\}$ and define the following \emph{discrete dynamical system}: For $x_{0} \in \mathbb{R}^{n}$, define
\begin{equation} 
\label{matdyn}
x_{j + 1}\ \ :=\ \ A_{\omega_{j}} x_{j}, \quad \text{for}\ \omega_{j} \in \big\{ 1,\, 2,\, \cdots,\, N \big\}. 
\end{equation}
That is, from a point $x_{j}$ at time $t = j$, we arrive at the point $x_{j + 1}$ at time $t = j + 1$ in the iteration of any generic point in $\mathbb{R}^{n}$, by randomly choosing one of the matrices from the above mentioned finite collection and the action by the chosen matrix. Observe that in order to achieve proper meaning to the above mentioned iterative scheme, one expects to understand nonhomogeneous products of matrices. 
\medskip 

\noindent 
Recall that given a self map $f$ on a topological space $X$, an element $x \in X$ is called a \emph{periodic point} of $f$ if there exists a positive integer $q$ such that $f^{q} (x) = x$. In such a case, the smallest such integer $q$ that satisfies $f^{q} (x) = x$ is called the \emph{period} of the periodic point $x$. The starting point of this work is the following consequence of the Perron-Frobenius theorem.
\medskip 

\noindent 
\begin{theorem} 
\label{starting-theorem}
Let $A \in M_{n} (\mathbb{R}_{+})$. Then, there exists a positive integer $q$ such that for every $x \in \mathbb{R}^{n}$ with $\left( \Vert A^{k} x \Vert \right)_{k\, \in\, \mathbb{N}}$ bounded, we have 
\[ \lim\limits_{k\, \to\, \infty} A^{k q} x\ \ =\ \ \xi_{x}, \] 
where $\xi_{x}$ is a periodic point of $A$ whose period divides $q$. 
\end{theorem}
\medskip 

\noindent 
We are interested in a generalization of Theorem \eqref{starting-theorem}, when the matrix $A$ in the above theorem is replaced by a product of the matrices $A_{r}$'s, possibly an infinite one, drawn from the finite collection of nonnegative matrices, $\big\{ A_{1},\, \ldots,\, A_{N} \big\}$. Besides generalizing Theorem \eqref{starting-theorem} as described above, we also bring out the existence of common periodic points for the said collection of matrices.
\medskip 

\noindent 
This manuscript is organized as follows: In Section \eqref{main}, we introduce basic notations, however only as much necessary to state the main results of this paper, namely Theorems \eqref{result1nc}, \eqref{result2} and \eqref{q2}. In Section \eqref{ceivec}, we familiarise the readers with some results from the literature, on adequate conditions to impose on a collection of matrices that ensures the existence of common eigenvectors. In section \eqref{three}, we state and prove special cases of Theorem \eqref{result1nc} by assuming furthermore properties on the collection of matrices. These are written as Theorems \eqref{result1} and \eqref{result1nodia}. We follow this with a few examples, in the same section. We begin Section \eqref{noncommmat} with a few examples, that we construct using the ideas propounded in the proof of Theorem \eqref{result1}, prove Theorem \eqref{result1nc} and restate the same in an alternate setting. This is written as Theorem \eqref{result1ncalt}. In Section \eqref{five}, we prove Theorem \eqref{result2} and study the examples already encountered. In section \eqref{six}, we consider words of infinite length based on a finite collection of pairwise commuting matrices and write the proof of Theorem \eqref{q2}. 

\section{Main results} 
\label{main}

In this section, we introduce some notations and explain the underlying settings of the main results and state our main results of this paper. As explained in the introductory section, we fix a finite set of nonnegative matrices $\big\{ A_{1},\, \ldots,\, A_{N} \big\},\ N < \infty$. For any finite $M \in \mathbb{N}$ and $p \in \mathbb{N}$, we denote the set of all $p$-lettered words on the set of first $M$ positive integers by 
\[ \Sigma_{M}^{p}\ \ :=\ \ \left\{ \omega = (\omega_{1} \omega_{2} \cdots \omega_{p})\ :\ \omega_{r} \in \big\{ 1,\, \cdots,\, M \big\} \right\}. \] 
For any $p$-lettered word $\omega := (\omega_{1} \omega_{2} \cdots \omega_{p}) \in \Sigma_{N}^{p}$, we define the (finite) matrix product 
\begin{equation} 
\label{aomega} 
A_{\omega}\ \ :=\ \ A_{\omega_{p}} \times A_{\omega_{p - 1}} \times \cdots \times A_{\omega_{2}} \times A_{\omega_{1}}. 
\end{equation}  

\noindent 
A key hypothesis in our first theorem assumes the existence of a nontrivial set of common eigenvectors, say $E' = \big\{ v_{1},\, v_{2},\, \ldots,\, v_{d} \big\}_{d\, \le\, n}$ for the given collection of matrices, $\big\{ A_{1},\, \ldots,\, A_{N} \big\}$. These common eigenvectors may be vectors in $\mathbb{R}^{n}$ or $\mathbb{C}^{n}$. As we will see in the next section, a sufficient condition that ensures the existence of common eigenvectors for the given collection is to demand the collection to be partially commuting, quasi-commuting or a Laffey pair when $N = 2$, or the collection to be quasi-commuting when $N \ge 3$. Define 
\[ \mathcal{LC} (E')\ \ =\ \ \big\{ \alpha_{1} v_{1}\, +\, \cdots\, +\, \alpha_{d} v_{d}\ :\ \alpha_{j} \in \mathbb{C}\ \text{satisfying}\ \alpha_{s_{1}} = \overline{\alpha_{s_{2}}}\ \forall v_{s_{1}} = \overline{v_{s_{2}}}\ \ \text{and}\ \ \alpha_{j} \in \mathbb{R}\ \text{otherwise} \big\}. \]
We now state our first result in this article. 
\medskip 

\begin{theorem} 
\label{result1nc}
Let $\big\{ A_{1},\, A_{2},\, \cdots,\, A_{N} \big\},\ N < \infty$, be a collection of $n \times n$ matrices with nonnegative entries, each having spectral radius $1$. Assume that the collection satisfies at least one of the following conditions, that ensures the existence of a nontrivial set of common eigenvectors. 
\begin{enumerate} 
\item If $N = 2$, the collection is either partially commuting, quasi-commuting or a Laffey pair. 
\item If $N \ge 3$, the collection is quasi-commuting. 
\end{enumerate} 
Let $E'$ denote the set of all common eigenvectors of the collection of matrices. For any finite $p$, let $\omega \in \Sigma_{N}^{p}$ and $A_{\omega}$ be the matrix associated to the word $\omega$. Then, for any vector $x \in \mathcal{LC} (E')$, there exists an integer $q_{\omega} \geq 1$ such that 
\begin{equation} 
\label{Limit}
\lim\limits_{k\, \to\, \infty} A_{\omega}^{k q_{\omega}} x\ \ =\ \ \xi_{(x,\, \omega)}, 
\end{equation}
where $\xi_{(x,\, \omega)}$ is a periodic point of $A_{\omega}$, whose period divides $q_{\omega}$. Moreover, when $p \geq N$ and $\omega$ is such that for all $1 \leq r \leq N$, there exists $1 \leq j \leq p$ such that $\omega_{j} = r$, the integer $q_{\omega}$ and the limiting point $\xi_{(x,\, \omega)} \in \mathbb{R}^{n}$ are independent of the choice of $\omega$.
\end{theorem}  
\medskip 

\noindent 
We now denote the interior of the nonnegative orthant of $\mathbb{R}^{n}$ by $(\mathbb{R}^{n}_{+})^{\circ}$, a convex cone and define the logarithm map and the exponential map, that appear frequently in nonlinear Perron-Frobenius theory as follows: $\log : (\mathbb{R}^{n}_{+})^{\circ} \longrightarrow \mathbb{R}^{n}$ and $\exp : \mathbb{R}^{n} \longrightarrow (\mathbb{R}^{n}_{+})^{\circ}$ by  
\begin{equation} 
\log (x)\ \ =\ \ \left( \log x_{1},\, \cdots,\, \log x_{n} \right)\ \ \ \ \text{and}\ \ \ \ \exp (x)\ \ =\ \ \left( e^{x_{1}},\, \cdots,\, e^{x_{n}} \right). 
\end{equation} 

\noindent 
As one may expect, these functions act as inverses of each other in the interior of $\mathbb{R}^{n}_{+}$. More on these functions and their uses in nonlinear Perron-Frobenius theory can be found in the monograph \cite{LN-book}. A nonnegative matrix, when viewed as a linear map on $\mathbb{R}^{n}$, preserves the partial order induced by $\mathbb{R}^{n}_{+}$. A map $f$ defined on a cone in $\mathbb{R}^{n}$ is said to be \emph{subhomogeneous} if for every $\lambda \in [0, 1]$, we have $\lambda f(x) \le f( \lambda x)$ for every $x$ in the cone and \emph{homogeneous} if $f ( \lambda x) = \lambda f(x)$ for every nonnegative $\lambda$ and every $x$ in the cone. It is then easy to verify that the function $f\ :=\ \exp \circ A \circ \log$ is a well-defined subhomogeneous map on $(\mathbb{R}^{n}_{+})^{\circ}$. Assuming the hypotheses of Theorem \eqref{result1nc}, we now state the second theorem of this paper for an appropriate subhomogeneous map, $f_{\omega}$. 
\medskip 

\begin{theorem} 
\label{result2} 
Let $\big\{ A_{1},\, \cdots,\, A_{N} \big\},\ N < \infty$ be a set of $n \times n$ matrices satisfying all the hypotheses in Theorem \eqref{result1nc}. For any finite $p$, let $\omega \in \Sigma_{N}^{p}$ and $A_{\omega}$ be the matrix associated with the word $\omega$. Consider the function $f_{\omega} : (\mathbb{R}^{n}_{+})^{\circ} \longrightarrow (\mathbb{R}^{n}_{+})^{\circ}$ by $f_{\omega} = \exp \circ A_{\omega} \circ \log$. Then, for any $y = e^{x} \in (\mathbb{R}^{n}_{+})^{\circ}$ where $x \in \mathcal{LC} (E')$, there exists an integer $q \geq 1$ such that 
\begin{equation} 
\label{limit}
\lim\limits_{k\, \to\, \infty} f_{\omega}^{k q} y\ \ =\ \ \eta_{y}, 
\end{equation}
where $\eta_{y}$ is a periodic point of $f_{\omega}$, whose period divides $q$. 
\end{theorem}
\medskip 

\noindent 
Our final theorem in this paper concerns the orbit of some $x \in \mathbb{R}^{n}$ under the action of some infinitely long word, whose letters belong to $\big\{ A_{1},\, \ldots,\, A_{N} \big\}$. In order to make our lives simpler, we shall assume that the given collection of matrices are pairwise commuting, with each matrix being diagonalizable over $\mathbb{C}$. This ensures the existence of $n$ linearly independent common eigenvectors $E' = \big\{ v_{1},\, \cdots,\, v_{n} \big\}$ for the given collection. Let the first $\kappa$ of these common eigenvectors correspond to eigenvalues of modulus $1$ for every matrix $A_{r}$, in the collection. 
\medskip 

\noindent 
For any $p$-lettered word $\omega \in \Sigma_{N}^{p}$ that has the presence of all $N$ letters, we denote by $\overline{\omega}$, the infinite-lettered word obtained by concatenating $\omega$ with itself, infinitely many times, {\it i.e.}, $\overline{\omega} = \left( \omega\, \omega\, \ldots \right)$. We know, from Theorem \eqref{result1nc}, that upon satisfying the necessary technical conditions, $\lim\limits_{k\, \to\, \infty} A_{\omega}^{k q} x = \xi_{x}$. Thus, the following definition makes sense. Let 
\[ \widetilde{A_{\omega}}\ \ :=\ \ A_{\overline{\omega}}\ \ :=\ \ \left( A_{\omega}^{q} \right)^{k},\ \text{as}\ k \to \infty. \] 
However, since Theorem \eqref{result1nc} only asserts $\xi_{x}$ to be a periodic point whose period divides $q$, we shall define $\widetilde{A_{\omega}} : \mathbb{R}^{n} \longrightarrow \left( \mathbb{R}^{n} \right)^{q}$. The precise action of $\widetilde{A_{\omega}}$ on points in $\mathbb{R}^{n}$ is given by 
\begin{equation} 
\label{aomegatilde} 
\widetilde{A_{\omega}} (x)\ \ :=\ \ \left( \xi_{x},\, A_{\omega} \xi_{x},\, \cdots,\, A_{\omega}^{q - 1} \xi_{x} \right). 
\end{equation} 

\noindent 
Let $\tau = \big( \tau_{1}\, \tau_{2}\, \tau_{3}\, \cdots \big) \in \big\{ 1,\, \ldots,\, N \big\}^{\mathbb{N}} =: \Sigma_{N}$ be any arbitrary infinite lettered word that encounters all the $N$ letters within a finite time, say $m$. It is easy to observe that the sequence 
\[ \Big( \overline{( \tau_{1}\, \cdots\, \tau_{m} )},\; \overline{( \tau_{1}\, \cdots\, \tau_{m + 1} )},\;\cdots \Big)\ \ \ \ \text{converges to $\tau$},  \] 
when $\Sigma_{N}$ is equipped with the usual product metric. For any $p \ge m$, denote by $\overline{\tau^{[p]}}$, the infinite-lettered word $\overline{( \tau_{1}\, \cdots\, \tau_{p} )}$ that occurs in the sequence, as described above that converges to any given $\tau$. Moreover, from the discussion above, we have that 
\[ \widetilde{A_{\tau^{[p]}}} x\ \ =\ \ \left( \xi_{x},\; A_{\tau^{[p]}} \xi_{x},\; \cdots,\; A_{\tau^{[p]}}^{q - 1} \xi_{x} \right). \] 
Notice that the first component of the vector in $\left( \mathbb{R}^{n} \right)^{q}$ is always $\xi{_x}\ \forall p \ge m$. Further, we define for every $r \in \big\{ 1,\, \cdots,\, N \big\}$, 
\[ \Phi_{(\tau,\, r)} (p)\ \ =\ \ \#\ \text{of}\ A_{r}\ \text{in}\ A_{\tau^{[p]}}. \] 

\noindent 
We now state our third theorem in this paper. 
\medskip 

\noindent 
\begin{theorem} 
\label{q2} 
Let $\big\{ A_{1},\, \cdots,\, A_{N} \big\},\ N < \infty$, be a collection of $n \times n$ pairwise commuting matrices with nonnegative entries, each having spectral radius $1$ and each matrix being diagonalizable over $\mathbb{C}$. Suppose $\tau \in \Sigma_{N}$ encounters all the $N$ letters within a finite time, say $m$. Then, for any $x \in \mathbb{R}^{n}$, there exists an increasing sequence $\big\{ p_{\gamma} \big\}_{\gamma\, \ge\, 1}$ of positive integers and a finite collection of positive integers $\big\{ \Lambda_{(r, j)} \big\}$ for $ 1 \le r \le N$ and $1 \le j \le \kappa$ such that 
\[ \sum_{r\, =\, 1}^{N} \Lambda_{(r, j)} \left[ \Phi_{(\tau,\, r)} (p_{\gamma_{k}}) - \Phi_{(\tau,\, r)} (p_{\gamma_{k'}}) \right]\ \ \equiv\ \ 0 \pmod q,\ \ \ \ \forall 1 \le j \le \kappa, \] 
where $p_{\gamma_{k}}$ and $p_{\gamma_{k'}}$ are any two integers from the sequence $\big\{ p_{\gamma} \big\}$. 
\end{theorem} 

\section{Common eigenvectors for a collection of matrices} 
\label{ceivec}

\noindent 
A key ingredient in one of our main results in this work is the existence of a nontrivial set of common eigenvectors for a given collection $\big\{ A_{1},\, \ldots,\, A_{N} \big\}$ of matrices. It is a well known result that if every matrix in the collection is diagonalizable over $\mathbb{C}$ with the collection commuting pairwise, there is a common similarity matrix that puts all the matrices in a diagonal form. A collection of non-commuting matrices may or may not have common eigenvectors. The question as to which collections of matrices possess common eigenvectors is extremely nontrivial. In what follows, we give a brief account of this question that is essential for this work. We begin with the following definition.
\medskip 

\noindent 
\begin{definition} 
\label{qc}
A collection $\big\{ A_{1},\, \ldots,\, A_{N} \big\}$ of matrices is said to be \emph{quasi-commuting} if for each pair $(r,\, s)$ of indices, both $A_{r}$ and $A_{s}$ commute with their (additive) commutator $[A_{r},\, A_{s}] := A_{r} A_{s} - A_{s} A_{r}$. 
\end{definition}
\medskip 

\noindent 
A classical result of McCoy (Theorem (2.4.8.7), \cite{HJ-book}) says the following. 
\medskip 

\begin{theorem} 
\label{mccoy}
Let $\big\{ A_{1},\, \ldots,\, A_{N} \big\}$ be a collection of $n \times n$ matrices. The following statements are equivalent. 
\begin{enumerate}
\item For every polynomial $p\, (t_{1},\, \ldots,\, t_{N})$ in $N$ non-commuting variables $t_{1},\, \ldots,\, t_{N}$ and every $r, s = 1,\, \ldots,\, N,\ p\, (A_{1},\, \ldots,\, A_{N}) [A_{r},\, A_{s}]$ is nilpotent.
\item There is a unitary matrix $U$ such that $U^{\ast} A_{r} U$ is upper triangular for every $r = 1, \ldots, N$.
\item There is an ordering $\lambda_{1}^{(r)}, \ldots, \lambda_{n}^{(r)}$ of the eigenvalues of each of the matrices $A_{r},\ 1 \le r \le N$ such that for any polynomial $p\, (t_{1},\, \ldots,\, t_{N})$ in $N$ non-commuting variables, the eigenvalues of $p\, (A_{1},\, \ldots,\, A_{N})$ are $p\, \left( \lambda_{j}^{(1)},\, \ldots,\, \lambda_{j}^{(N)} \right),\ j = 1, \ldots, n$.
\end{enumerate}
\end{theorem} 
\medskip 

\noindent 
If the matrices and the polynomials are over the real field, then all calculations may be carried out over $\mathbb{R}$, provided all the matrices have eigenvalues in $\mathbb{R}$. It turns out that a sufficient condition that guarantees any of the above three statements is when the collection of matrices is quasi-commutative (see Drazin {\it et al} \cite{drazin-dungey-gruenberg}). Moreover, the first statement implies that the collection $\big\{ A_{1},\, \ldots,\, A_{N} \big\}$ has common eigenvectors. There are also other classes of matrices which possess common eigenvectors. 
\medskip 

\noindent 
A pair $(A,\, B)$ of matrices is said to \emph{partially commute} if they have common eigenvectors. Moreover, two matrices $A$ and $B$ partially commute {\it iff} the \emph{Shemesh subspace} $\mathcal{N} = \displaystyle \bigcap_{k, l\, =\, 1}^{n - 1} {\rm ker} \left( \left[ A^{k},\, B^{l} \right] \right)$ is a nontrivial maximal invariant subspace of $A$ and $B$ over which both $A$ and $B$ commute (see Shemesh, \cite{Sh}). A pair $(A,\, B)$ of matrices is called a \emph{Laffey pair} if ${\rm rank} \left( [A,\, B] \right) = 1$. It can be shown that such a pair of matrices partially commute, but do not commute.

\section{A special case of Theorem \eqref{result1nc}}
\label{three} 

\noindent 
In this section, we first state and prove a special case of Theorem \eqref{result1nc} by assuming furthermore properties for the collection of matrices, $\big\{ A_{1},\, \ldots,\, A_{N} \big\}$. Later, we observe that a proper modification of the proof of this special case yields a complete proof of Theorem \eqref{result1nc}, in the general case. 
\medskip 

\noindent 
We begin with the following result due to Frobenius. Suppose $A$ is an irreducible matrix in $M_{n} (\mathbb{R}_{+})$ such that there are exactly $\kappa$ eigenvalues of modulus $\rho (A)$. This integer $\kappa$ is called the \emph{index of imprimitivity} of $A$. If $\kappa = 1$, the matrix $A$ is said to be \emph{primitive}. If $\kappa > 1$, the matrix is said to be \emph{imprimitive}. The following result is due to Frobenius.
\medskip 

\noindent 
\begin{theorem}[\cite{Z-book}, Theorem (6.18)] 
\label{spectrum-nonnegative-matrix}
Let $A$ be an irreducible nonnegative matrix with its index of imprimitivity equal to $\kappa$. If $\lambda_{1},\, \ldots,\, \lambda_{\kappa}$ are the eigenvalues of $A$ of modulus $\rho(A)$, then $\lambda_{1},\, \ldots,\, \lambda_{\kappa}$ are the distinct $\kappa$-th roots of $[\rho (A)]^{\kappa}$. 
\end{theorem}
\medskip 

\noindent 
We now restate Theorem \eqref{result1nc} for the special case when the considered collection of matrices is pairwise commuting with each matrix being diagonalizable over $\mathbb{C}$ and prove the same. 
\medskip 

\noindent 
\begin{theorem} 
\label{result1} 
Let $\big\{ A_{1},\, \cdots,\, A_{N} \big\},\ N < \infty$, be a collection of $n \times n$ pairwise commuting matrices with nonnegative entries, each having spectral radius $1$ and each matrix being diagonalizable over $\mathbb{C}$. For any finite $p$, let $\omega \in \Sigma_{N}^{p}$ and $A_{\omega}$ be the matrix associated to the word $\omega$. Then, for any $x \in \mathbb{R}^{n}$, there exists an integer $q_{\omega} \geq 1$ such that 
\begin{equation} 
\label{Limit}
\lim\limits_{k\, \to\, \infty} A_{\omega}^{k q_{\omega}} x\ \ =\ \ \xi_{(x,\, \omega)}, 
\end{equation}
where $\xi_{(x,\, \omega)}$ is a periodic point of $A_{\omega}$, whose period divides $q_{\omega}$. Moreover, when $p \geq N$ and $\omega$ is such that for all $1 \leq r \leq N$, there exists $1 \leq j \leq p$ such that $\omega_{j} = r$, the integer $q_{\omega}$ and the limiting point $\xi_{(x,\, \omega)} \in \mathbb{R}^{n}$ are independent of the choice of $\omega$.
\end{theorem}
\medskip 

\noindent 
\begin{proof} 
We first observe that the hypotheses in the statement of the Theorem ensures that the matrices $A_{1},\, \ldots,\, A_{N}$ are simultaneously diagonalizable. Let $E' = \big\{ v_{1},\, \ldots,\, v_{n} \big\}$ be a set of $n$ linearly independent common eigenvectors of the matrices $A_{1},\, \ldots,\, A_{N}$ that satisfies $A_{r} v_{s}\ =\ \lambda_{(r,\, s)} v_{s}$, where $\lambda_{(r,\, s)}$ is an eigenvalue of the matrix $A_{r}$ corresponding to the eigenvector $v_{s}$. Observe that for any $p$-lettered word $\omega = \left( \omega_{1} \cdots \omega_{p} \right)$, we have  
\[ A_{\omega} v_{s}\ \ = \ \ \lambda_{(\omega_{p},\, s)} \cdots \lambda_{(\omega_{1},\, s)} v_{s}\ \ = \ \ \lambda_{(\omega,\, s)} v_{s},\ \ \text{where}\ \ \lambda_{(\omega,\, s)}\ = \ \lambda_{(\omega_{p},\, s)} \cdots \lambda_{(\omega_{1},\, s)}. \] 

\noindent 
We now rearrange the common eigenvectors $\big\{ v_{1},\, \ldots,\, v_{n} \big\}$ as $\big\{ v_{1},\, \ldots,\, v_{\kappa},\, v_{\kappa + 1},\, \ldots,\, v_{n} \big\}$, where $\kappa$ is defined as 
\begin{equation} 
\label{kappa} 
\kappa\ \ =\ \ \# \left\{ v_{s}\ :\ A_{r} v_{s}\ =\ \lambda_{(r,\, s)} v_{s}\ \text{with}\ \left\vert \lambda_{(r,\, s)} \right\vert = 1\ \forall 1 \leq r \leq N \right\}.
\end{equation}

\noindent 
Let $S_{2}$ be a subset of $\Sigma_{\kappa}^{2}$, defined as 
\begin{eqnarray*} 
S_{2} & := & \big\{ s = ( s_{1} s_{2} ) \in \Sigma_{n}^{2}\ :\ v_{s_{1}} = \overline{v_{s_{2}}} \big\}\ \ \ \text{and} \\ 
\mathcal{LC} (E') & := & \big\{ \alpha_{1} v_{1} + \cdots + \alpha_{n} v_{n}\ :\ \alpha_{j} \in \mathbb{C}\ \text{satisfying}\ \alpha_{s_{1}} = \overline{\alpha_{s_{2}}}\ \forall s \in S_{2}\ \ \text{and}\ \ \alpha_{j} \in \mathbb{R}\ \text{otherwise} \big\}. 
\end{eqnarray*} 
Note by our definition that $\mathcal{LC} (E') \subset \mathbb{R}^{n}$. Owing to the hypotheses on the spectral radius in the statement of the theorem, we have that for every $x \in \mathbb{R}^{n}$, the sequence $\big\{ \| A_{\omega}^{k} x \| \big\}_{k\, \geq\, 1}$ is bounded. In fact, 
\[ \big\Vert A_{\omega}^{k} x \big\Vert\ \ =\ \ \big\Vert \alpha_{1} A_{\omega}^{k} v_{1}\; +\; \cdots\; +\; \alpha_{n} A_{\omega}^{k} v_{n} \big\Vert\ \ \leq\ \ \big\vert \alpha_{1} \big\vert \big\Vert v_{1} \big\Vert\; +\; \cdots\; +\; \big\vert \alpha_{n} \big\vert \big\Vert v_{n} \big\Vert. \] 

\noindent 
Let $q_{1},\, \ldots,\, q_{N}$ be positive integers that satisfies the outcome of Theorem \eqref{starting-theorem}, for the matrices $A_{1},\, \ldots,\, A_{N}$ respectively. For some $p > N$, let $\omega$ be a $p$-lettered word in $\Sigma_{N}^{p}$ such that for all $1 \leq r \leq N$, there exists $1 \leq j \leq p$ such that $\omega_{j} = r$. Define $q$ to be the least common multiple of the numbers $\lbrace q_{1},\, \ldots,\, q_{N} \rbrace$. 
\medskip 

\noindent 
For every $s \in \big\{ 1,\, \ldots,\, n \big\}$ and $r \in \big\{ 1,\, \ldots,\, N \big\}$, we enumerate the following possibilities that can occur for the values of $\lambda_{(r,\, s)}$: 
\begin{enumerate}
\item[{\bf Case 1.}] $\left( \lambda_{(r,\, s)} \right)^{q} = 1$ for every $r$ and for some $s$ with $\lambda_{(r,\, s)} \in \mathbb{R}$. This implies that the corresponding eigenvector $v_{s}$ lies in $\mathcal{LC} (E')$. 
\item[{\bf Case 2.}] $\left( \lambda_{(r,\, s)} \right)^{q} = 1$ for every $r$ and for some $s$ with $\lambda_{(r,\, s)} \in \mathbb{C}$. This implies that there exists eigenvectors $v_{s}$ and $\overline{v_{s}}$ with corresponding eigenvalues conjugate to each other such that $\alpha_{s} v_{s} + \overline{\alpha_{s} v_{s}}$ lies in $\mathcal{LC} (E')$. 
\item[{\bf Case 3.}] $\big\vert \lambda_{(r,\, s)} \big\vert < 1$ for some $s$ and for some $r$. In this case, the iterates of $v_{s}$ under the map $A_{\omega}$ goes to $0$; that is, $\lim\limits_{k\, \to\, \infty} A_{\omega}^{k} v_{s} = 0$. 
\end{enumerate}

\noindent 
For any $x \in \mathbb{R}^{n}$ that can be written as $x\, =\, \alpha_{1} v_{1}\, +\, \cdots\, +\, \alpha_{n} v_{n}$, we have  
\begin{eqnarray*} 
\lim\limits_{k\, \to\, \infty} A_{\omega}^{k q} x & = & \alpha_{1} \lim\limits_{k\, \to\, \infty} \left( \lambda_{(\omega,\, 1)} \right)^{k q} v_{1}\; +\; \cdots\; +\; \alpha_{n} \lim\limits_{k\, \to\, \infty} \left( \lambda_{(\omega,\, n)} \right)^{k q} v_{n} \\ 
& = & \alpha_{1} v_{1}\; +\; \cdots\; +\; \alpha_{\kappa} v_{\kappa} \\ 
& =: & \xi_{(x,\, \omega)}. 
\end{eqnarray*}

\noindent 
Observe that $\xi_{(x,\, \omega)}$ and $q$ are independent of the length of the word and in fact, the word $\omega$ itself. We denote $\xi_{(x,\, \omega)} = \xi_{x}$. Moreover, $\xi_{x} \in \mathcal{LC} (E')$. 
\medskip 

\noindent 
Further, for any point $v = \beta_{1} v_{1}\, +\, \ldots\, +\, \beta_{\kappa} v_{\kappa} \in \mathcal{LC} (E')$ and $1 \le r \le N$, we have 
\[ A_{r}^{q_{r}} v\ \ =\ \ \beta_{1} A_{r}^{q_{r}} v_{1}\; +\; \cdots\; +\; \beta_{\kappa} A_{r}^{q_{r}} v_{\kappa}\ \ =\ \ \beta_{1} \left( \lambda_{(r,\, 1)} \right)^{q_{r}} v_{1}\; +\; \cdots\; +\; \beta_{\kappa} \left( \lambda_{(r,\, \kappa)} \right)^{q_{r}} v_{\kappa}\ \ =\ \ v. \] 
Since $\xi_{x} \in \mathcal{LC} (E')$, it is a periodic point of $A_{1}, \ldots, A_{N}$ with periods $q_{1}, \ldots, q_{N}$ respectively. Thus, $\xi_{x}$ is a periodic point of $A_{\omega}$ with period $q$. 
\end{proof}
 \bigskip 

\noindent 
We now explore the possibility of weakening the diagonalizability condition in the hypothesis of Theorem \eqref{result1}. However, since the matrices in the collection $\big\{ A_{1},\, \ldots,\, A_{N} \big\}$ commute pairwise, we still obtain a collection of common eigenvectors, say $E' = \big\{ v_{1},\, \ldots,\, v_{d} \big\}$. For example, consider the following pair of non-diagonalizable, commuting matrices: 
\[ A = \begin{bmatrix} 1 & 1 \\ 0 & 1 \end{bmatrix}\ \ \ \ \text{and}\ \ \ \ B = \begin{bmatrix} 1 & 2 \\ 0 & 1 \end{bmatrix}. \] 
Observe that $e_{1} = (1, 0)^{t}$ is a common eigenvector for $A$ and $B$, whereas $e_{2} = (0, 1)^{t}$ is a common generalized eigenvector for $A$ and $B$. We further note that for any word $\omega$ that contains both the letters, the orbit of $e_{2}$ under $A_{\omega}$ is unbounded while that of $e_{1}$ is bounded. 
\medskip 

\noindent 
Thus, in order to demand the orbit of a point under some relevant map $A_{\omega}$ to be bounded, we concern ourselves only with the common eigenvectors, say $E'$ for the collection of matrices $\big\{ A_{1},\, \ldots,\, A_{N} \big\}$. Suppose $d$ is the cardinality of $E'$ while $\kappa$ is the cardinality of those eigenvectors in $E'$ with corresponding eigenvalues being of modulus $1$ for every matrix in the collection. Then, the argument in the proof of Theorem \eqref{result1} goes through verbatim for any $x \in \mathcal{LC} (E')$. Observe that for any $\omega \in \Sigma_{N}^{p}$, we have the orbit of any such $x$ under $A_{\omega}$, {\it i.e.}, $\displaystyle{\left\{ \Vert A_{\omega}^{k} x \Vert \right\}_{k\, \ge\, 1}}$ to be bounded. Thus, we have proved: 
\medskip 

\begin{theorem} 
\label{result1nodia}
Let $\big\{ A_{1},\, \ldots,\, A_{N} \big\},\ N < \infty$, be a collection of $n \times n$ pairwise commuting matrices with nonnegative entries, each having spectral radius $1$. For any finite $p$, let $\omega \in \Sigma_{N}^{p}$ and $A_{\omega}$ be the matrix associated to the word $\omega$. Then, for any $x \in \mathcal{LC} (E')$, the same conclusion as in Theorem \eqref{result1} holds. 
\end{theorem} 
\medskip 

\noindent 
In other words, if $x$ is an element of $\mathbb{R}^{n}$ whose orbit under $A_{\omega}$ is norm bounded, then $\displaystyle \lim_{k\, \to\, \infty} A_{\omega}^{k} x$ exists and is a periodic point of $A_{\omega}$, with its period dividing the order of a permutation on $N$ symbols. In the absence of diagonalizability, the boundedness of the orbit of an element of $\mathbb{R}^{n}$ can be guaranteed only on a nontrivial subset of $\mathbb{R}^{n}$. Moreover, observe from the proof that $\xi_{x}$ is a common periodic point for all the matrices in the collection. 
\medskip 

\noindent 
\begin{remark}
A few remarks are in order.	 
\begin{enumerate}
\item If all the matrices $A_{1},\, \ldots,\, A_{N}$ are pairwise commuting symmetric matrices, subject to the spectral radius assumption in Theorem \eqref{result1}, then the periods of all the periodic points corresponding to the eigenvalues $1$ and $- 1$ for all $A_{r}$'s is $2$. Hence, for a matrix product corresponding to a word $\omega$, we have $q = 2$. 
\item We have proved Theorems \eqref{result1} and \eqref{result1nodia} for a special choice of $\omega$ that contains all the $N$ letters. Suppose $\omega'$ is any arbitrary $p$-lettered word. Then we can take the appropriate subset of $\big\{ 1,\, \ldots,\, N \big\}$, whose members have been used for the writing of the word $\omega'$ and the same result as above follows for $\omega'$. 
\end{enumerate}
\end{remark}
\medskip 

\noindent 
Suppose the $p$-lettered word $\omega = \left(r r \cdots r \right)$ for some $1 \leq r \leq N$. Then the above theorem reduces to a particular case of a result of Lemmens as stated in \cite{L}. We now state the same as a corollary.
\medskip 

\noindent 
\begin{corollary} 
\label{lemmens-cor-1}
Let $A$ be an $n \times n$ matrix with nonnegative entries that is diagonalizable over $\mathbb{C}$ and of spectral radius $1$. Then there exists an integer $q \ge 1$ such that for every $x \in \mathbb{R}^{n}$, we have $\lim\limits_{k\, \to\, \infty} A^{q k} x = \xi_{x}$, where $\xi_{x}$ is a periodic point of $A$ with its period dividing $q$.
\end{corollary}
\medskip 

\noindent 
We now present a few examples that illustrate Theorem \eqref{result1}. We fix a few notations before proceeding further with the examples. We denote the standard basis vectors of $\mathbb{R}^{n}$ by $e_{1},\, \ldots,\, e_{n}$, while $I_{n}$ denotes the identity matrix of order $n$. We write the permutation matrices in column partitioned form; for instance, we denote the $2 \times 2$ permutation matrix $\big[\, e_{2}\; |\; e_{1}\, \big]$ by $J_{2}$. The matrix of $1$'s (of any order) is denoted by $J$. The diagonal matrix of order $n$ with diagonal entries $d_{1},\, \ldots,\, d_{n}$ is denoted by ${\rm diag} (d_{1},\, \cdots,\, d_{n})$. Our first example is a fairly simple one and illustrates the scenario in Corollary \eqref{lemmens-cor-1}. 
\medskip 

\noindent 
\begin{example} 
\label{example-1}
Consider the diagonalizable matrix $A = J_{2}$ with spectral radius $1$. If $x = e_{2} \in \mathbb{R}^{2}$, then observe that 
\[ A^{k} x\ \ =\ \ \begin{cases} e_{2}, & \text{if $k$ is even}, \\ e_{1}, & \text{if $k$ is odd}. \end{cases} \] 
In this example, we obtain $q = 2$.
\end{example}
\medskip 

\noindent 
The second example involves a pair of $6 \times 6$ commuting nonnegative matrices.
\medskip 

\noindent 
\begin{example} 
\label{example-2}
Consider 
\begin{displaymath} 
\begin{array}{r c l c r c l c r c l} 
A & = & \begin{bmatrix} A_{1} & 0 \\ 0 & A_{2} \end{bmatrix} & \text{where} & A_{1} & = & \big[\, e_{4}\; |\; e_{1}\; |\; e_{2}\; |\; e_{3}\, \big] & \text{and} & A_{2} & = & \displaystyle{\frac{1}{3}} \begin{bmatrix} 1 & 2 \\ 2 & 1 \end{bmatrix}; \\ 
B & = & \begin{bmatrix} B_{1} & 0 \\ 0 & B_{2} \end{bmatrix} & \text{where} & B_{1} & = & \displaystyle{\frac{1}{4}} \begin{bmatrix} 0 & 2 & 0 & 2 \\ 2 & 0 & 2 & 0 \\ 0 & 2 & 0 & 2 \\ 2 & 0 & 2 & 0  \end{bmatrix} & \text{and} & B_{2} & = & \displaystyle{\frac{1}{10}} \begin{bmatrix} 3 & \sqrt{7} \\ \sqrt{7} & 3 \end{bmatrix}. 
\end{array} 
\end{displaymath} 

\noindent 
It can be easily seen that the matrices $A$ and $B$ commute and are diagonalizable over $\complex$ and therefore, are simultaneously diagonalizable. Further, 
\[ {\rm spec}(A)\ \ =\ \ \left\{ 1,\, 1,\, - 1,\, i,\, - i,\, - \frac{1}{3} \right\}\ \ \ \text{and}\ \ \ {\rm spec}(B)\ \ =\ \ \left\{ 0,\, 0,\, 1,\, - 1,\, \frac{3 + \sqrt{7}}{10},\, \frac{3 - \sqrt{7}}{10} \right\}. \]
The following table gives the common eigenvectors of $A$ and $B$ and the corresponding eigenvalues of the matrices $A$ and $B$. 
\[ \begin{tabular}{|| c || c | c | c | c | c | c ||}
\hline \hline 
\text{Eigenvectors} & $v_{1}$ & $v_{2}$ & $v_{3}$ & $v_{4}$ & $v_{5}$ & $v_{6}$ \\ 
\hline 
\text{Eigenvalues of $A$} & $1$ & $- 1$ & $- i$ & $i$ & $1$ & $- 1/3$ \\ 
\hline 
Eigenvalues of $B$ & $1$ & $- 1$ & $0$ & $0$ & $ \lambda_{1}$ or $\lambda_{2}$ & $\lambda_{2}$ or $\lambda_{1}$ \\ 
\hline \hline  
\end{tabular} \]
where $\lambda_{j} = \displaystyle{\frac{3 \pm \sqrt{7}}{10}},\ j = 1, 2$. Moreover, the first four eigenvectors are given by 
\begin{displaymath} 
\begin{array}{r c l r c l} 
v_{1} & = & (1,\, 1,\, 1,\, 1,\, 0,\, 0)^{t}, & v_{2} & = & (1,\, - 1,\, 1,\, - 1,\, 0,\, 0)^{t}, \\ 
v_{3} & = & (1,\, i,\, - 1,\, - i,\, 0,\, 0)^{t}, & v_{4} & = & (1,\, - i,\, - 1,\, i,\, 0,\, 0)^{t}. 
\end{array} 
\end{displaymath}  
As we shall see below, the eigenvectors $v_{5}$ and $v_{6}$ do not play any role in our analysis and hence, we do not write the same here. Let $x \in \mathbb{R}^{6}$. If $x = \alpha_{1} v_{1}$ with $\alpha_{1} \in \mathbb{R}$, then $x$ is a fixed point for $A$ and $B$, and hence is a fixed point for $A^{p_{1}} B^{p_{2}}$ for any $p_{1},\, p_{2} \in \mathbb{N}$. If $x = \alpha_{2} v_{2}$ with $\alpha_{2} \in \mathbb{R}$, then $x$ has period $2$ for both $A$ and $B$, and hence is a fixed point for $A^{p_{1}} B^{p_{2}}$ for any $p_{1},\, p_{2} \in \mathbb{N}$ that satisfies $p_{1} + p_{2} \equiv 0 \pmod 2$. If $x = \alpha_{3} v_{3} + \alpha_{4} v_{4}$ with $\alpha_{3}, \alpha_{4} \in \mathbb{C}$ satisfying $\alpha_{3} = \overline{\alpha_{4}}$, then $x$ has period $4$ for $A$, whereas $B x = 0$. If $x = \alpha_{5} v_{5} + \alpha_{6} v_{6}$ with $\alpha_{5}, \alpha_{6} \in \mathbb{R}$, then $\displaystyle \lim_{k\, \to\, \infty} A^{k} x = \alpha_{5} v_{5}$ and $\displaystyle \lim_{k\, \to\, \infty} B^{k} x = 0$. Hence, for any $p_{1},\, p_{2} \in \mathbb{N}$, we have $\displaystyle \lim_{k\, \to\, \infty} \left( A^{p_{1}} B^{p_{2}} \right)^{k} x = 0$. 
\medskip 

\noindent
Thus, for any word $A_{\omega}$ that contains both $A$ and $B$, the periodic point $x = \alpha_{1} v_{1} + \alpha_{2} v_{2}$ has period $2$, whereas the least common multiple of the periods of the periodic points of $A$ and $B$ is $4$.
\end{example}
\medskip 

\noindent 
We now present another example, this time in $\mathbb{R}^{7}$, illustrating our result.
\medskip 

\noindent 
\begin{example} 
\label{example-3}
Let 
\begin{displaymath} 
\begin{array}{r c l c r c l c r c l} 
A & = & \begin{bmatrix} I_{3} & 0 & 0 \\ 0 & J_{2} & 0 \\ 0 & 0 & D_{A}\end{bmatrix} & \text{where} & D_{A} & = & {\rm diag} \left( \displaystyle{\frac{1}{2},\, \frac{1}{3}} \right); & & & & \\ 
\\ 
B & = & \begin{bmatrix} J_{3} & 0 & 0 \\ 0 & I_{2} & 0 \\ 0 & 0 & D_{B} \end{bmatrix} & \text{where} & D_{B} & = & {\rm diag} \left( \displaystyle{\frac{1}{5},\, \frac{1}{6}} \right) & \text{and} & J_{3} & = & \big[\, e_{3}\; |\; e_{1}\; |\; e_{2}\, \big]. 
\end{array} 
\end{displaymath} 

\noindent 
Clearly, $A B = B A$ and $A$ and $B$ are diagonalizable and hence, are simultaneously diagonalizable. As earlier,  we write a table with the common eigenvectors and the corresponding eigenvalues for the matrices $A$ and $B$. 
\[ \begin{tabular}{|| c || c | c | c | c | c | c | c ||}
\hline \hline  
\text{Eigenvectors} & $v_{1}$ & $v_{2}$ & $v_{3}$ & $v_{4}$ & $v_{5}$ & $v_{6}$ & $v_{7}$\\ 
\hline 
\text{Eigenvalues of $A$} & $1$ & $1$ & $1$ & $1$ & $-1$ & $- 1/2$ & $- 1/3$ \\ 
\hline 
Eigenvalues of $B$ & $1$ & $\omega$ & $\omega^{2}$ & $1$ & $1$ & $1/5$ & $1/6$\\ 
\hline \hline  
\end{tabular} \] 
where $\omega$ is the cubic root of unity and the $v_{i}$'s are 
\begin{displaymath} 
\begin{array}{r c l r c l} 
v_{1} & = & (1,\, 1,\, 1,\, 0,\, 0,\, 0,\, 0)^{t} & v_{2} & = & (1,\, \omega,\, \omega^{2},\,  0,\, 0,\, 0,\, 0)^{t} \\ 
v_{3} & = & (1,\, \omega^{2},\, \omega,\,  0,\, 0,\, 0,\, 0)^{t} & v_{4} & = & (0,\, 0,\, 0,\, 1,\, 1,\, 0,\, 0)^{t} \\ 
v_{5} & = & (0,\, 0,\, 0,\, 1,\, - 1,\, 0,\, 0)^{t} & v_{6} & = & e_{6},\ \ \ \ v_{7}\ \ =\ \ e_{7}. 
\end{array} 
\end{displaymath} 
Consider $x \in \mathbb{R}^{7}$. If $x = \alpha_{1} v_{1} + \alpha_{4} v_{4}$ with $\alpha_{1}, \alpha_{4} \in \mathbb{R}$, then $A x = x = B x$. If $x = \alpha_{2} v_{2} + \alpha_{3} v_{3}$ with $\alpha_{2}, \alpha_{3} \in \mathbb{C}$ satisfying $\alpha_{2} = \overline{\alpha_{3}}$, then $A x = x$ and $B^{3} x = x$. If $x = \alpha_{5} v_{5}$, then $A^{2} x = x$ and $B x = x$. If $x = \alpha_{6} v_{6} + \alpha_{7} v_{7}$, then $\displaystyle \lim_{k\, \to\, \infty} A^{k} x = 0$ and $\displaystyle \lim_{k\, \to\, \infty} B^{k} x = 0$. If $x = \alpha_{2} v_{2} + \alpha_{3} v_{3} + \alpha_{5} v_{5}$, with $\alpha_{2} = \overline{\alpha_{3}}$ and $\alpha_{5} \in \mathbb{R}$, then $A^{2} x = x,\; B^{3} x = x,\; (A B)^{6} x = x$ and $\left( A^{p_{1}} B^{p_{2}} \right)^{6} x = x$ for every $p_{1},\, p_{2} \in \mathbb{N}$ with $p_{1} \equiv 1 \pmod 2$ and $p_{2} \equiv 1 \pmod 3$ or $2 \pmod 3$. 
\medskip

\noindent
In this case, for any word $A_{\omega}$ that contains both $A$ and $B$, the periodic point $x = \alpha_{1} v_{1} + \alpha_{2} v_{2} + \alpha_{3} v_{3} + \alpha_{4} v_{4} + \alpha_{5} v_{5}$, with $\alpha_{2}, \alpha_{3} \in \mathbb{C}$ satisfying $\alpha_{2} = \overline{\alpha_{3}}$ and $\alpha_{i} \in \mathbb{R}$ for $i = 1, 4, 5$ has period $6$, which is the least common multiple of the periods of the periodic points of $A$ and $B$. 
\end{example}

\section{What happens when the matrices do not commute?} 
\label{noncommmat}

\noindent 
We begin this section with a few examples in the non-commuting set up that satisfy the hypotheses in Theorem \eqref{result1nc}. We observe from the proof of Theorem \eqref{result1} that the common eigenvectors of the matrices concerned play a significant role in our analysis. Thus, we concern ourselves with matrices that have a nontrivial set of common eigenvectors, even when they do not commute. 
\medskip 

\noindent 
\begin{example}
\label{example-4}
Let 
\begin{displaymath} 
\begin{array}{r c l c r c l c r c l} 
A & = & \begin{bmatrix} J_{4} & 0 \\ 0 & A' \end{bmatrix} & \text{where} & A' & = & \begin{bmatrix} 1/5 & 1/6 \\ 1/6 & 1/5 \end{bmatrix} & \text{and} & J_{4} & = & \big[\, e_{4}\; |\; e_{1}\; |\; e_{2}\; |\; e_{3}\, \big]; \\ 
\\
B & = & \begin{bmatrix} J_{2} & 0 & 0 \\ 0 & J_{2} & 0 \\ 0 & 0 & B' \end{bmatrix} & \text{where} & B' & = & \begin{bmatrix} 1/7 & 1/8 \\ 1/7 & 1/8 \end{bmatrix}. & & & & 
\end{array} 
\end{displaymath} 

\noindent 
Here, $A B \neq B A$. The common eigenvectors of the matrices $A$ and $B$ are given by 
\begin{displaymath} 
\begin{array}{r c l r c l} 
v_{1} & = & \left( 1,\, 1,\, 1,\, 1,\, 0,\, 0 \right)^{t}, & v_{2} & = & \left( 1,\, - 1,\, 1,\, - 1,\, 0,\, 0 \right)^{t},  \\  
v_{3} & = & \left( 0,\, 0,\, 0,\, 0,\, 1,\, 1 \right)^{t}, & v_{4} & = & \left( 1,\, i,\, - 1,\, - i,\, 0,\, 0 \right)^{t} \\ 
v_{5} & = & \left( 1,\, - i,\, - 1,\, i,\, 0,\, 0 \right)^{t}. & & & 
\end{array} 
\end{displaymath} 

The corresponding eigenvalues of the respective matrices are given in the following table. 
\[ \begin{tabular}{|| c || c | c | c | c | c ||}
\hline \hline 
\text{Eigenvectors} & $v_{1}$ & $v_{2}$ & $v_{3}$ & $v_{4}$ & $v_{5}$ \\ 
\hline 
\text{Eigenvalues of $A$} & $1$ & $- 1$ & $1/5 + 1/6$ & $- i$ & $i$ \\ 
\hline 
Eigenvalues of $B$ & $1$ & $- 1$ & $1/7 + 1/8$ & $0$ & $0$ \\ 
\hline \hline  
\end{tabular} \]

\noindent 
The non-common eigenvector of $A$ is given by $v_{A} = \left( 0,\, 0,\, 0,\, 0,\, 1,\, - 1 \right)^{t}$ with the corresponding eigenvalue being $(1/5 - 1/6)$, while the non-common eigenvector of $B$ is given by $v_{B} = \left( 0,\, 0,\, 0,\, 0,\, 7,\, - 8 \right)^{t}$ with the corresponding eigenvalue being $0$. Suppose $x \in \mathbb{R}^{6}$. Then, $A^{k} v_{A}$ and $B^{k} v_{B}$ converge to $0$ as $k \to \infty$. Further, if $x = \alpha_{1} v_{1}$ with $\alpha_{1} \in \mathbb{R}$, then $x$ is a fixed point of $A$ and $B$. If $x = \alpha_{2} v_{2}$ with $\alpha_{2} \in \mathbb{R}$, then $x$ is a periodic point of both $A$ and $B$ with period $2$. If $x = \alpha_{3} v_{3}$ with $\alpha_{3} \in \mathbb{R}$, then $\displaystyle \lim_{k\, \to\, \infty} A^{k} x = 0$ and $\displaystyle \lim_{k\, \to\, \infty} B^{k} x = 0$. If $x = \alpha_{4} v_{4} + \alpha_{5} v_{5}$ with $\alpha_{4}, \alpha_{5} \in \mathbb{C}$ satisfying $\alpha_{4} = \overline{\alpha_{5}}$, then $x$ is a periodic point of $A$ with period $4$, whereas $B x = 0$.
\medskip

\noindent 
Thus, for any word $A_{\omega}$ that contains both $A$ and $B$, the vector $x = \alpha_{1} v_{1} + \alpha_{2} v_{2}$ is a periodic point with period $2$, whereas the least common multiple of the periods of the periodic points is $4$.
\end{example}  
\medskip 

\noindent 
\begin{example}
\label{example-5}
Let 
\[ A\ =\ \begin{bmatrix} I_{3} & 0 & 0 \\ 0 & J_{2} & 0 \\ 0 & 0 & \frac{1}{2} J \end{bmatrix}\ \ \text{and}\ \ B\ =\ \begin{bmatrix} J_{3} & 0 & 0 \\ 0 & I_{2} & 0 \\ 0 & 0 & B' \end{bmatrix}\ \ \text{where}\ \ B'\ =\ \begin{bmatrix} 1/3 & 1/4 \\ 1/3 & 1/4 \end{bmatrix}. \] 

\noindent 
In this example, $A$ and $B$ do not commute; however, they have the following six common eigenvectors: 
\[ \begin{array}{r c l r c l r c l} 
v_{1} & = & \left( 1,\, \omega,\, \omega^{2},\, 0,\, 0,\, 0,\, 0 \right)^{t}, & v_{2} & = & \left( 1,\, \omega^{2},\, \omega,\, 0,\, 0,\, 0,\, 0 \right)^{t}, & v_{3} & = & \left( 1,\, 1,\, 1,\, 0,\, 0,\, 0,\, 0 \right)^{t}, \\ 
v_{4} & = & \left( 0,\, 0,\, 0,\, 1,\, 1,\, 0,\, 0 \right)^{t}, & v_{5} & = & \left( 0,\, 0,\, 0,\, 1,\, - 1,\, 0,\, 0 \right)^{t}, & v_{6} & = & \left(0,\, 0,\, 0,\, 0,\, 0,\, 1,\, 1 \right)^{t}. 
\end{array} \] 

\noindent 
As earlier, we write the corresponding the eigenvalues of the matrices in the following table: 
\[ \begin{tabular}{|| c || c | c | c | c | c | c ||}
\hline \hline 
\text{Eigenvectors} & $v_{1}$ & $v_{2}$ & $v_{3}$ & $v_{4}$ & $v_{5}$ & $v_{6}$ \\ 
\hline 
\text{Eigenvalues of $A$} & $1$ & $1$ & $1$ & $1$ & $- 1$ & $1$ \\ 
\hline 
Eigenvalues of $B$ & $\omega$ & $\omega^{2}$ & $1$ & $1$ & $1$ & $1/3 + 1/4$ \\ 
\hline \hline  
\end{tabular} \]

\noindent
In this case, for any word $A_{\omega}$ that contains both $A$ and $B$, the vector $x = \alpha_{1} v_{1} + \alpha_{2} v_{2} + \alpha_{3} v_{3} + \alpha_{4} v_{4} + \alpha_{5} v_{5}$ with $\alpha_{1}, \alpha_{2} \in \mathbb{C}$ satisfying $\alpha_{1} = \overline{\alpha_{2}}$ and $\alpha_{3}, \alpha_{4}, \alpha_{5} \in \mathbb{R}$ is a periodic point with period $6$, which is the least common multiple of the periods of the periodic points of $A$ and $B$.
\end{example} 
\medskip 

\noindent 
At this juncture, we write two further examples that illustrate different scenarios in the non-commuting set-up, that are worth observing. The first one focuses on an eigenvector of the product $AB$, which is not a common eigenvector for $A$ and $B$. We draw the attention of the readers to the fact that the common eigenvectors of the matrices $A$ and $B$, in the above mentioned examples, have a bounded orbit under the action of $A_{\omega}$ that contains both $A$ and $B$. In the first example, we show that this does not happen for a specific vector. In the second example, we observe that the non-common eigenvectors of a pair of non-commuting matrices can have a bounded orbit. However, the corresponding periodic points may depend on the word, even if both the matrices are present in the word. 
\medskip 

\noindent 
\begin{example} 
\label{example-6} 
Consider 
\[ A\ \ =\ \ \frac{1}{3} \begin{bmatrix} 3 & 0 & 0 \\ 0 & 1 & 2 \\ 0 & 2 & 1 \end{bmatrix}\ \ \ \ \text{and}\ \ \ \ B\ \ =\ \ \frac{1}{3} \begin{bmatrix} 3 & 0 & 0 \\ 0 & 1 & 4 \\ 0 & 1 & 1 \end{bmatrix}. \] 

\noindent 
Observe the matrices $A$ and $B$ do not commute, are both diagonalizable, have spectral radius $1$ and with one common eigenvector. However, the spectral radius of the product $A B$ is equal to $(2 + \sqrt{3})/3$. Suppose $v$ is the eigenvector corresponding to the above mentioned eigenvalue of $A B$. Then, observe that the sequence $\left( \| (AB)^{k} v \| \right)_{k\, \in\, \mathbb{N}}$ is unbounded. Hence, $\displaystyle{\lim_{k\, \to\, \infty} (AB)^k v}$ does not exist for the word $AB$. 
\end{example} 
\medskip

\noindent
\begin{example}\label{example-7}
Let 
\[ A\ \ =\ \ \frac{1}{3} \begin{bmatrix} 1 & 2 \\ 2 & 1 \end{bmatrix}\ \ \ \ \text{and}\ \ \ \ B\ \ =\ \ \frac{1}{5} \begin{bmatrix} 1 & 4 \\ 2 & 3 \end{bmatrix}. \] 
It can be easily seen that $AB \ne BA$. The eigenvalues of $A$ and $B$ are $1, - \frac{1}{3}$ and $1, - \frac{1}{5}$ respectively. The vector $(1, 1)^{t}$ is a common eigenvector for $A$ and $B$ corresponding to the eigenvalue $1$. Moreover, the eigenvalues of $A B$ are $1, \frac{1}{15}$ (and so the same is true for $B A$). It easily follows from this that any $x \in \mathbb{R}^{2}$ has a bounded orbit. The eigenvector corresponding to the eigenvalue $- \frac{1}{3}$ for $A$ is $(1, -1)^{t}$ and the eigenvector corresponding to the eigenvalue $- \frac{1}{5}$ for $B$ is $(2, -1)^{t}$. Note that 
\[ (A B)^{k}\ -\ (B A)^{k}\ \ =\ \ \begin{bmatrix} - \alpha(k) & \alpha(k) \\ -\alpha(k) & \alpha(k) \end{bmatrix}\ \ \ \ \text{for}\ \ k \ge 1, \] 
and therefore the commutator has rank $1$, making this a Laffey pair. It is now obvious that 
\[ \lim_{k\, \to\, \infty} (A B)^{k} (1, 1)^{t}\ \ =\ \ \lim_{k\, \to\, \infty} (B A)^{k} (1, 1)^{t}\ \ \ \ \text{since}\ \ \ \ ((A B)^{k} - (B A)^{k}) (1, 1)^{t}\ \ =\ \ (0, 0)^{t}. \] 
Nevertheless, 
\begin{eqnarray*} 
((A B)^{k} - (B A)^{k}) (2, -1)^{t} & = & - 3 \alpha(k) (1, 1)^{t}\ \ \ \ \text{whereas} \\ 
((A B)^{k} - (B A)^{k}) (1, -1)^{t} & = & - 2 \alpha(k) (1, 1)^{t}. 
\end{eqnarray*} 
Therefore, if $x$ is one of the points $(2,-1)^t$ or $(1,-1)^t$, then, $\displaystyle \lim_{k \to \infty} (AB)^k x \neq \displaystyle \lim_{k \to \infty} (BA)^k x$, since $\displaystyle{\lim_{k\, \to\, \infty} \alpha(k) \ne 0}$. 
\end{example}
\medskip 

\noindent 
\begin{proof}[of Theorem \eqref{result1nc}] By our hypothesis, there exists a nontrivial set of common eigenvectors for the given collection of matrices, $\big\{ A_{1},\, \ldots,\, A_{N} \big\}$, say $E' = \big\{ v_{1},\, \ldots,\, v_{d} \big\}$. Analogous to the definitions of $\kappa,\, S_{2}$ and $\mathcal{LC} (E')$ in the proof of Theorem \eqref{result1}, we now define: 
\begin{eqnarray*} 
\kappa & := & \# \left\{ v_{s} \in V\ :\ A_{r} v_{s}\ =\ \lambda_{(r,\, s)} v_{s}\ \text{with}\ \left\vert \lambda_{(r,\, s)} \right\vert = 1\ \forall 1 \leq r \leq N \right\} \\ 
S_{2} & := & \big\{ s = ( s_{1} s_{2} ) \in \Sigma_{d}^{2}\ :\ v_{s_{1}} = \overline{v_{s_{2}}} \big\}\ \ \ \text{and} \\ 
\mathcal{LC} (E') & := & \big\{ \alpha_{1} v_{1} + \cdots + \alpha_{d} v_{d}\ :\ \alpha_{j} \in \mathbb{C}\ \text{satisfying}\ \alpha_{s_{1}} = \overline{\alpha_{s_{2}}}\ \forall s \in S_{2}\ \ \text{and}\ \ \alpha_{j} \in \mathbb{R}\ \text{otherwise} \big\}. 
\end{eqnarray*} 
A summary of our findings from the above-stated examples illustrate that norm-boundedness of orbits of points in $\mathbb{R}^{n}$ holds only on $\mathcal{LC} (E')$. Hence, the remainder of the proof follows along the same lines as in the proof of Theorem \eqref{result1}. 
\end{proof} 
\medskip 

\noindent 
We now describe another way of writing Theorem \eqref{result1nc}. Recall that $\Sigma_{N}$ denotes the set of all infinite-lettered words on the set of symbols $\big\{ 1,\, \ldots,\, N \big\}$. Considering the Cartesian product of the symbolic space $\Sigma_{N}$ and $\mathbb{R}^{n}$, one may describe the dynamical system discussed in this paper thus: Given a collection $\big\{ A_{1},\, \ldots,\, A_{N} \big\}$ of $n \times n$ matrices, let $T : X = \Sigma_{N} \times \mathbb{R}^{n} \longrightarrow X$ be defined by $T \left( \left( \tau,\, x \right) \right) = \left( \sigma \tau,\, A_{\tau_{1}} x \right)$ where $\tau = \big( \tau_{1} \tau_{2} \tau_{3} \ldots \big)$ and $\sigma$ is the shift map defined on $\Sigma_{N}$ by $\left( \sigma \tau \right)_{n} = \tau_{n + 1}$ for $n \ge 1$. The discrete topology accorded on the set of symbols defines a product topology on $\Sigma_{N}$, thereby making $\Sigma_{N}$ a compact and perfect metric space. We then equip $X$ with the corresponding topology and study $T$ as a non-invertible map. 
\medskip 

\noindent 
\begin{theorem} 
\label{result1ncalt}
Let $\big\{ A_{1},\, \cdots,\, A_{N} \big\},\ N < \infty$, be a collection of $n \times n$ matrices with nonnegative entries, each having spectral radius $1$. Assume that the collection satisfies at least one of the following conditions, that ensures the existence of a nontrivial set of common eigenvectors. 
\begin{enumerate} 
\item If $N = 2$, the collection is either partially commuting, quasi-commuting or a Laffey pair. 
\item If $N \ge 3$, the collection is quasi-commuting. 
\end{enumerate} 
Let $E'$ denote the set of all common eigenvectors of the collection of matrices. Let $\tau \in \Sigma_{N}$ be any arbitrary infinite lettered word that encounters all the $N$ letters within a finite time, say $m$. For $p \ge m$, let $\left\{ \overline{\tau^{[p]}} \right\}$ be a sequence of infinite-lettered words that converges to $\tau$. Let $A_{\tau^{[p]}}$ be the matrix associated to the $p$-lettered word $\tau^{[p]} \in \Sigma_{N}^{p}$. Then, for every $p \ge m$ and any vector $x \in \mathcal{LC} (E')$, there exists an integer $q \geq 1$ such that 
\begin{equation} 
\lim\limits_{k\, \to\, \infty} T^{k p q} ( \tau,\, x )\ \ =\ \ ( \tau,\, \xi_{x} ), 
\end{equation}
where $( \tau,\, \xi_{x} )$ is a periodic point of $T$, whose period divides the least common multiple of $p$ and $q$. 
\end{theorem}

\section{Proof of Theorem \eqref{result2}} 
\label{five} 

\noindent 
In this section, we prove Theorem \eqref{result2} and illustrate the result for the examples mentioned in the earlier sections. Recall the definitions of the logarithm map and the exponential map that were defined component-wise from Section \eqref{main}. If $x$ in an element in the interior of $\mathbb{R}^{n}_{+}$, then it can be easily verified that 
\[ f(x)\ \ =\ \ \exp \circ A \circ \log (x)\ \ =\ \ \left( \prod_{j\, =\, 1}^{n}  x_{j}^{a_{1 j} },\; \prod_{j\, =\, 1}^{n} x_{j}^{a_{2 j}},\; \cdots,\; \prod_{j\, =\, 1}^{n}  x_{j}^{a_{n j}}\ \right). \] 

\noindent 
We  point out that a generalization of Theorem \eqref{starting-theorem} holds for continuous order-preserving subhomogeneous maps on a polyhedral cone (Theorem 8.1.7, \cite{LN-book}). A specific example of such a map is the above defined function $f = \exp \circ A \circ \log$, where $A$ is a nonnegative matrix. 
\medskip 

\noindent 
If $A$ is a nonnegative matrix, then $f$ can be extended continuously to a homogeneous map on $\mathbb{R}^{n}_{+}$. Since each entry $a_{ij}$ is nonnegative, it follows that for all $1 \leq i, j \leq n,\ a_{ij}\neq 0$ implies $x_{j}^{a_{ij}} \to 0$ whenever $x_{j} \to 0$;  if $a_{ij} = 0$ then $x_{j}^{a_{ij}} = 1$. Therefore, the above expression has limits on the boundary of the nonnegative orthant. It can also be verified that when $A$ is a row stochastic matrix, the function $f$ is a homogeneous map on $\mathbb{R}^{n}_{+}$.
\medskip

\noindent
We now write the proof of Theorem \eqref{result2}. 
\medskip 

\noindent 

\noindent 
\begin{proof}[of Theorem \eqref{result2}] 
The proof of this theorem is similar to that of Theorem \eqref{result1}. Let $f_{1},\, \cdots,\, f_{N}$ be self-maps on $(\mathbb{R}^{n}_{+})^{\circ}$ given by $f_{r}\ =\ \exp\circ A_{r} \circ \log\ \forall 1 \le r \le N$. By Theorem \eqref{starting-theorem}, for every $A_{r}$ there exists $q_{r}$ such that for every $x \in \mathbb{R}^{n}$ with a corresponding norm-bounded orbit, we have $\lim\limits_{k\, \to\, \infty} A_{r}^{k q_{r}} x = \xi_{(x, r)}$, where $\xi_{(x, r)}$ is a periodic point of $A_{r}$. Consider $y  = e^{x} \in (\mathbb{R}^{n}_{+})^{\circ}$. Then, 
\[ \lim\limits_{k\, \to\, \infty} f_{r}^{k q_{r}} y\ \ =\ \ \lim\limits_{k\, \to\, \infty} \exp \circ A_{r}^{k q_{r}} \circ \log \left( e^{x} \right) \ \ =\ \ \exp \left( \lim\limits_{k\, \to\, \infty} A_{r}^{k q_{r}} x \right)\ \ =\ \ \exp(\xi_{(x, r)}). \]

\noindent 
It is clear that $\exp(\xi_{(x, r)})$ is a periodic point of $f_{r}$, since 
\[ f_{r}^{q_{r}} (\exp(\xi_{(x, r)}))\ \ =\ \ ( \exp \circ A_{r} \circ \log )^{q_{r}} ( \exp ( \xi_{(x, r)}))\ \ =\ \ \exp \circ A_{r}^{q_{r}} (( \xi_{(x, r)}))\ \ =\ \ \exp ( \xi_{(x, r)}). \]

\noindent 
As earlier, we define $q$ to be the least common multiple of $\big\{ q_{1},\, q_{2},\, \cdots,\, q_{N} \big\}$. Then, for a given $p$-lettered word $\omega$ and for every $y = e^{x} \in (\mathbb{R}^{n}_{+})^{\circ}$ where $x \in \mathcal{LC} (E')$, 
\[ \lim\limits_{k\, \to\, \infty} f_{\omega}^{k q} y\ \ =\ \ \lim\limits_{k\, \to\, \infty} \exp \circ A_{\omega}^{k q} \circ \log \left( e^{x} \right)\ \ =\ \ \exp \left( \lim\limits_{k\, \to\, \infty} A_{\omega}^{q k} x \right)\ \ =\ \ \exp( \xi_{(x, \omega)})\ \ =:\ \ \eta_{(y, \omega)}. \]

\noindent 
One can prove that  $\eta_{(y, \omega)}$ is a periodic point of $f_{\omega}$ exactly along the same lines as the proof of $\eta_{(y, r)} = \exp(\xi_{(x, r)})$ being a periodic point of $f_{r}$. Further, the independence of $\eta_{(y, \omega)}$ from the $p$-lettered word $\omega$ can be established, as in the proof of Theorem \eqref{result1}. 
\end{proof} 
\medskip 

\noindent 
\begin{remark} 
In addition to the hypotheses of Theorem \eqref{result2}, suppose each of the matrices in the collection $\big\{ A_{1},\, \ldots,\, A_{N} \big\}$ is row stochastic ({\it i.e.}, the row sum being $1$ for every row). Let $x$ be a point on the diagonal of $(\mathbb{R}_{n}^{+})^{\circ}$; that is, $x = \left( x_{1},\, \ldots,\, x_{n} \right)$ where $x_{1} = \cdots = x_{n} > 0$. Then, an easy calculation yields $\xi_{x} = \eta_{x}$. 
\end{remark} 
\medskip 

\noindent 
\begin{example} 
\label{example-8}
For the matrices given in Example \eqref{example-2}, it is easily seen that $f_{A}$ and $f_{B}$ are 
\begin{eqnarray*}
f_{A} (x) & = & \left( x_{2},x_{3},x_{4},x_{1},\sqrt[3]{x_{5}x_{6}^{2}},\sqrt[3]{x_{5}^{2}x_{6}}\right)\quad \text{and } \\ 
f_{B} (x) & = & \left(\sqrt{x_{2}x_{4}},\sqrt{x_{1}x_{3}},\sqrt{x_{2}x_{4}},\sqrt{x_{1}x_{3}} ,\sqrt[3]{x_{5}^{2} x_{6}},\sqrt[3]{x_{5} x_{6}^{2}}\right)
\end{eqnarray*} 
Observe that the vector $\exp{(\alpha_{1} v_{1} + \alpha_{2} v_{2})} \in (\mathbb{R}^{6}_{+})^{\circ}$ with $\alpha_{i} \in \mathbb{R}$ is a periodic point with period $2$ for $f_{\omega}$, whenever $\omega$ is any word that contains both $A$ and $B$. 
\end{example}
\medskip 

\noindent 
\begin{example} 
\label{example-9}
For the matrices given in Example \eqref{example-3}, it is easily seen that $f_{A}$ and $f_{B}$ are 
\begin{equation*}
f_{A}\ =\ \left( x_{1},x_{2},x_{3},x_{5},x_{4},\sqrt{x_{6}},\sqrt[3]{x_{7}}\right)\ \ \ \ \text{and}\ \ \ \ f_{B}\ =\ \left( x_{2},x_{3},x_{1},x_{4},x_{5},\sqrt[5]{x_{6}},\sqrt[6]{x_{7}}\right). 
\end{equation*}
In this example, we note that the vector $\exp{(\alpha_{1} v_{1} + \alpha_{2} v_{2} + \alpha_{3} v_{3} + \alpha_{4} v_{4} + \alpha_{5} v_{5})} \in (\mathbb{R}^{7}_{+})^{\circ}$ with $\alpha_{2}, \alpha_{3} \in \mathbb{C}$ satisfying $\alpha_{2} = \overline{\alpha_{3}}$ and $\alpha_{1}, \alpha_{4}, \alpha_{5} \in \mathbb{R}$ is a periodic point for $f_{\omega}$ of period $6$, for any $\omega$ that contains both $A$ and $B$. 
\end{example}
\medskip 

\noindent
Now, we look at the examples discussed earlier in Section \eqref{noncommmat} in the non-commuting set up for the above defined subhomogeneous map. 
\medskip 

\noindent 
\begin{example}
\label{example-10} 
For the matrices given in Example \eqref{example-4}, we compute $f_{A}$ and $f_{B}$ as 
\begin{eqnarray*} 
f_{A} & = & \big( x_{2},\, x_{3},\, x_{4},\, x_{1},\, \sqrt[5]{x_{5}} \sqrt[6]{x_{6}},\, \sqrt[6]{x_{5}} \sqrt[5]{x_{6}} \big) \quad \text{and} \\ 
f_{B} & = & \big( \sqrt{x_{2} x_{4}},\, \sqrt{x_{1} x_{3}},\, \sqrt{x_{2} x_{4}},\, \sqrt{x_{1} x_{3}},\, \sqrt[7]{x_{5}} \sqrt[8]{x_{6}},\, \sqrt[7]{x_{5}} \sqrt[8]{x_{6}} \big). 
\end{eqnarray*} 
In this example, we note that the point $x = \exp{(\alpha_{1} v_{1} + \alpha_{2} v_{2})} \in (\mathbb{R}^{6}_{+})^{\circ}$ with $\alpha_{i} \in \mathbb{R}$ is a periodic point with period $2$ for $f_{\omega}$, where $\omega$ is any word that contains both the letters $A$ and $B$. 
\end{example}
\medskip 

\noindent 
\begin{example}
\label{example-11}
For the matrices given in Example \eqref{example-5}, we obtain $f_{A}$ and $f_{B}$ as 
\begin{eqnarray*} 
f_{A} & = & \big( x_{1},\, x_{2},\, x_{3},\, x_{5},\, x_{4},\, \sqrt{x_{6} x_{7}},\, \sqrt{x_{6} x_{7}} \big) \quad \text{and} \\ 
f_{B} & = & \big( x_{2},\, x_{3},\, x_{1},\, x_{4},\, x_{5},\, \sqrt[3]{x_{6}} \sqrt[4]{x_{7}},\, \sqrt[3]{x_{6}} \sqrt[4]{x_{7}} \big). 
\end{eqnarray*}
Here, the vector $\exp{(\alpha_{1} v_{1} + \alpha_{2} v_{2} + \alpha_{3} v_{3} + \alpha_{4} v_{4} + \alpha_{5} v_{5})} \in (\mathbb{R}^{7}_{+})^{\circ}$ with $\alpha_{1}, \alpha_{2} \in \mathbb{C}$ satisfying $\alpha_{1} = \overline{\alpha_{2}}$ and $\alpha_{3}, \alpha_{4}, \alpha_{5} \in \mathbb{C}$ is a periodic point for $f_{\omega}$ of period $6$, for any word $\omega$ that contains both $A$ and $B$. 
\end{example} 

\section{Words of infinite length} 
\label{six} 

\noindent 
In this section, we write the proof of Theorem \eqref{q2}. Recall that the hypotheses of Theorem \eqref{q2} and Theorem \eqref{result1} are one and the same. 
\medskip 

\noindent 
\begin{proof}[of Theorem \eqref{q2}] 
Recall from Equation \eqref{aomegatilde} that $\widetilde{A_{\tau^{[p]}}} x = \left( \xi_{x},\, A_{\tau^{[p]}} \xi_{x},\, \cdots,\, A_{\tau^{[p]}}^{q - 1} \xi_{x} \right)$ for $p \ge m$. Suppose $\widetilde{A_{\tau^{[p]}}} x = \widetilde{A_{\tau^{[p + k]}}} x\ \forall k \ge 0$, then consider any increasing subsequence of positive integers $\big\{ p_{\gamma} \big\}$, with $p_{\gamma} \ge m\ \forall \gamma \ge 1$. Observe that $\left\{ \widetilde{A_{\tau^{[p_{\gamma}]}}} x \right\}_{\gamma\, \ge\, 1}$ is a constant sequence in $\left( \mathbb{R}^{n} \right)^{q}$. 

\medskip 

\noindent 
In general, it is not necessary that $A_{\tau^{[m]}} \xi_{x} = A_{\tau^{[m + 1]}} \xi_{x}$. However, owing to $\xi_{x}$ being a periodic point of $A_{\tau^[p]}$ for $p \ge m$, whose period divides $q$ ($> 1$, say), a simple application of the pigeon-hole principle ensures $A_{\tau^{[m]}} \xi_{x} = A_{\tau^{[m']}} \xi_{x}$, for some $m' > m$. Moreover, pairwise commutativity of the matrices in the collection then, guarantees 
\[ \widetilde{A_{\tau^{[m]}}} x\ \ =\ \ \widetilde{A_{\tau^{[m']}}} x,\ \ \ \ \text{as vectors in}\ \left( \mathbb{R}^{n} \right)^{q}. \] 
Proceeding along similar lines, one obtains an increasing sequence, say $\big\{ p_{\gamma} \big\}$ such that $\left\{ \widetilde{A_{\tau^{[p_{\gamma}]}}} x \right\}_{\gamma\, \ge\, 1}$ is a constant sequence of vectors in $\left( \mathbb{R}^{n} \right)^{q}$. 
\medskip 

\noindent 
Choose any two integers $p_{\gamma_{k}}$ and $p_{\gamma_{k'}}$ from the sequence $\big\{ p_{\gamma} \big\}$. Then, $\widetilde{A_{\tau^{[p_{\gamma_{k}}]}}} x = \widetilde{A_{\tau^{[p_{\gamma_{k'}}]}}} x$. By a mere comparison of coordinates, we then obtain $A_{\tau^{[p_{\gamma_{k}}]}} \xi_{x} = A_{\tau^{[p_{\gamma_{k'}}]}} \xi_{x}$. Since $\xi_{x} = \displaystyle{\sum_{s\, =\, 1}^{\kappa} \alpha_{s} v_{s}}$, we obtain $\alpha_{1} \lambda_{(\tau^{[p_{\gamma_{k}}]}, 1)} v_{1}\, +\, \cdots\, +\, \alpha_{\kappa} \lambda_{(\tau^{[p_{\gamma_{k}}]}, \kappa)} v_{\kappa}\ =\ \alpha_{1} \lambda_{(\tau^{[p_{\gamma_{k'}}]}, 1)} v_{1}\, +\, \cdots\, +\, \alpha_{\kappa} \lambda_{(\tau^{[p_{\gamma_{k'}}]}, \kappa )} v_{\kappa}$. This implies that for every $1 \leq j \leq \kappa$, we have 
\begin{eqnarray*} 
\lambda_{(\tau^{[p_{\gamma_{k}}]}, j)}\ \ =\ \ \lambda_{(\tau^{[p_{\gamma_{k'}}]}, j)} & \Longleftrightarrow & \prod_{r\, =\, 1}^{N} \lambda_{(r, j)}^{\Phi_{(\tau, r)} (p_{\gamma_{k}})}\ \ =\ \ \prod_{r\, =\, 1}^{N} \lambda_{(r, j)}^{\Phi_{(\tau, r)} (p_{\gamma_{k'}})} \\ 
& \Longleftrightarrow & \prod_{r\, =\, 1}^{N} \lambda_{(r, j)}^{\Phi_{(\tau, r)} (p_{\gamma_{k}}) - \Phi_{(\tau, r)} (p_{\gamma_{k'}})}\ \ =\ \ 1. 
\end{eqnarray*} 

\noindent 
Since the numbers $\lambda_{(r, j)}$'s are $q^{\text{th}}$ roots of unity, we obtain positive integers $\Lambda_{(r, j)}$ that satisfies 
\[ \sum_{r\, =\, 1}^{N} \Lambda_{(r, j)} \left[ \Phi_{(\tau,\, r)} (p_{\gamma_{k}}) - \Phi_{(\tau,\, r)} (p_{\gamma_{k'}}) \right]\ \ \equiv\ \ 0 \pmod q,\ \ \ \ \forall 1 \le j \le \kappa. \] 
\end{proof}

\vspace{1cm}
\noindent

\bibliographystyle{amsplain}

\end{document}